\documentclass[12pt]{article}

\usepackage{amsfonts}

\newtheorem{thm}{Theorem}[section]

\newtheorem{cor}[thm]{Corollary}

\newtheorem{prop}[thm]{Proposition}

\newenvironment{remark}{\par\medskip\noindent{\bf Remark.\ }}{\par\smallskip}
\newcommand{\proof
}{\par\medskip\noindent {\bf Proof.\ \ }}

\newcommand{\be}{\begin{equation}}
\newcommand{\ee}{\end{equation}}
\newcommand{\openbox}{\leavevmode
  \hbox to8pt{\hfil\vrule\vbox to6pt{\hrule width6pt\vfil\hrule}\vrule}}

\newcommand{\ve}[1]{\mathbf{#1}}
\newcommand{\qed}{\hbox to5pt{ } \hfill \openbox\bigskip\medskip}

\newcommand{\Fp}{\mathbb F _p}
\newcommand{\Fq}{\mathbb F _q}

\newcommand{\N}{\mathbb N}

\newcommand{\F}{\mathbb F}

\title{Difference sets and  power residues}
\author{G\'abor Heged\H{u}s
\\{\normalsize  \'Obuda University}
\\{\normalsize B\'ecsi \'ut 96, Budapest, Hungary, H-1037}
\\{\normalsize hegedus.gabor@nik.uni-obuda.hu}
}

\begin{document}
\maketitle

\begin{abstract}
Let $p\geq 3$ be a prime and   $n\geq 1$ be  an integer.  Let $K\subseteq {\mathbb F _p}$ denote a fixed subset with $0\in K$. Let $A\subseteq ({\mathbb F _p})^n$ be an arbitrary subset such that 
$$
\{ \mathbf{a}-\mathbf{b}:~\mathbf{a},\mathbf{b}\in A,\mathbf{a}\neq \mathbf{b}\}\cap K^n=\emptyset.
$$
Then we prove the exponential upper bound
$$
|A|\leq ( p-|K|+ 1 )^n.
$$

We use in our proof the linear algebra bound method.
\end{abstract}
\medskip

\section{Introduction}


Let $p$ denote a prime with $p\equiv 1 \pmod 4$. The Paley graph  of order $p$ is a graph $G(p)$ on $p$ vertices (here we associate  each
vertex with an element of $\Fp$), where $(i, j)$ is an edge if $i-j$ is a quadratic residue modulo
$p$.  Let $\omega(p)$ denote the
clique number of the Paley graph of order $p$. It is a challenging open problem to determine $\omega(p)$.

Until now the best known upper bound  is $\omega(p)\leq \sqrt{p}-1$ for infinitely many primes $p$ (see \cite{BRM} Theorem 2.1).

It is well-known that the Paley graph is a self-complementary graph, hence $\alpha(G(p))=\omega(p)$. 

We can  consider the following reformulation of this problem:
let $Q(2)$ denote the set of quadratic residues in $\Fp$. How large can a set $A\subseteq \Fp$ be given that 
$$
\{a-b:~ a,b\in A,a\neq b\}\subseteq \Fp \setminus Q(2)?
$$

We investigate here the following generalization of this problem to elementary $p$-groups. 
Let $p\geq 3$ be a prime, $k\geq 2$ be a fixed integer and  let $Q(k)$ denote the set of $k$th power residues modulo $p$ (i.e. $Q(k)=\{b\in \Fp:~\mbox{ there exists }x\in \Fp\mbox{ with } x^k\equiv b \pmod p\}$. Clearly $0\in Q(k)$. 
Let $n\geq 1$ be a fixed integer. How large can a set $A\subseteq ({\Fp})^n$ be given that 
$$
\{\ve a-\ve b:~\ve a,\ve b\in A,\ve a\neq \ve b\}\subseteq ({\Fp})^n\setminus (Q(k))^n?
$$

Matolcsi and Ruzsa investigated  the following version of this question in \cite{MR}:

Let $G$ denote a finite Abelian group and let $B\subseteq G$ be a fixed standard set (i.e. $B=-B$ and $0\in B$). Consider the number 
$$
\Delta(B):=\max \{|A|:~ A\subseteq G, (A-A)\cap B=\{0\}\}.
$$
How large can $\Delta(B)$ be for a a fixed standard set?

Green investigated a similar question  in \cite{G}. For the reader's convenience we state here his result. 
Let $q$ be a prime power and $n\geq 1$ be a fixed integer. Denote by $P(q,n)$ the $n$-dimensional vector space over the finite field $\Fq$ of all polynomials $a_{n-1}x^{n-1}+\ldots +a_0$ of degree less than $n$.
Green proved  the following result in \cite{G}.

\begin{thm} \label{green}
Let $k\geq 2$ be a fixed integer and let $q$ be a prime power. Define
$$
c(k,q):=(2k^2D_q(k)^2\log(q))^{-1},
$$ 
where $D_q(k)$ is the sum of the digits of $k$ in base $q$. Suppose that $A\subseteq P(q,n)$ is a subset with $|A|> 2q^{(1-c(k,q))n}$, then $A$ contains distinct polynomials $p(x)$ and  $q(x)$ such that $p(x)-q(x)= h(x)^k$ for some $h(x)\in \Fq[x]$.
\end{thm}

We state here our main results.
 
\begin{thm} \label{main}
Let $p\geq 3$ be a prime and  let $n\geq 1$ be a fixed integer. Let $K\subseteq {\Fp}$ be a fixed subset with $0\in K$. Define $t:=|K|$. Suppose that  $A\subseteq ({\Fp})^n$ is a subset such that 
$$
|A|> (p-t + 1)^n. 
$$
Then there exist $\ve a_1,\ve a_2\in A$, $\ve a_1\neq \ve a_2$ such that $\ve a_1-\ve a_2 \in K^n$.
\end{thm}
\begin{remark}
We think the bound $(p-t + 1)^n$ is  not optimal in general.
The only obvious case, when our bound  is sharp, is the following:
Let $K:=\Fp$.  Then $t=p$ and clearly if $A\subseteq ({\Fp})^n$ is an arbitrary subset with $|A|>1$, then there exist $\ve a_1,\ve a_2\in A$, $\ve a_1\neq \ve a_2$ such that $\ve a_1-\ve a_2 \in K^n=(\Fp)^n$.
\end{remark}

Our proof technique is the usual linear algebra bound method (see \cite{BF}  Chapter 2). 
Finally we point out an important  special case of Theorem \ref{main}.

\begin{cor} \label{main2}
Let $p\geq 3$ be a prime, $k\geq 2$ be a fixed integer and let $Q(k)\subseteq {\Fp}$ denote the set of $k$th power residues modulo $p$. Let $n\geq 1$ be a fixed integer. Define $d:=gcd(k,p-1)$.   Suppose that $A\subseteq ({\Fp})^n$ is a subset such that 
\begin{equation} \label{dpower}
|A|> \Big( \frac{(p-1)(d-1)}{d} + 1 \Big)^n. 
\end{equation} 
Then there exist $\ve a_1,\ve a_2\in A$, $\ve a_1\neq \ve a_2$ such that $\ve a_1-\ve a_2 \in (Q(k))^n$.
\end{cor}
\proof

Define $K:=Q(k)$. 

 It is a well-known fact that $t:=|Q(k)|=\frac{p-1}{d}+1$ (see \cite{N} Theorem 3.11). Hence $p-t=p-(\frac{p-1}{d}+1)=(p-1)-\frac{p-1}{d}=\frac{(p-1)(d-1)}{d}$ and Corollary \ref{main2} follows from Theorem \ref{main}. \qed

\begin{remark}
Consider the special case $n=1$ in Theorem \ref{main}. Then it follows from Theorem  \ref{main} that if 
$$
|A|>2q^{(1-c(k,q))},
$$
where
$$
c(k,q):=(2k^2D_q(k)^2\log(q))^{-1},
$$ 
then there exist $a_1,a_2\in A$, $a_1\neq a_2$ such that $a_1-a_2 \in Q(k)$. This bound is clearly better than our bound appearing in the inequality  (\ref{dpower}), but it works only in the case $n=1$.
\end{remark}

\section{Proof}

We can prove our main result using the linear algebra bound method and the Determinant Criterion (see \cite{BF} Proposition 2.7). We recall here for the reader's convenience the Determinant Criterion.
\begin{prop} \label{det} (Determinant Criterion)
Let $\F$ denote an arbitrary field. Let $f_i:\Omega \to \F$ be functions for each $i=1,\ldots ,m$ and $\ve v_i\in \Omega$ elements such that the $m \times m$ matrix $B=(f_i(\ve v_j))_{i,j=1}^m$
is nonsingular. Then $f_1,\ldots ,f_m$ are linearly independent functions of the space $\F^{\Omega}$. 
\end{prop}

{\bf Proof of Theorem \ref{main}:}\\

Indirectly, suppose that there exists an $A\subseteq ({\Fp})^n$  subset such that 
$$
|A|> (p-t + 1)^n 
$$
and 
\begin{equation} \label{indir}
\{\ve a-\ve b:~\ve a,\ve b\in A,\ve a\neq \ve b\}\subseteq ({\Fp})^n\setminus K^n.
\end{equation}

Define $N:={\Fp}\setminus K$. Then $|N|=p-t.$

Consider the polynomial
$$
Q(x_1,\ldots ,x_n):=\prod_{1\leq i\leq n} \prod_{\alpha\in N} (x_i-\alpha)\in \Fp[x_1, \ldots , x_n].
$$
Then clearly 
$$
\mbox{deg}(Q)=n|N|=n(p-t).
$$

If we expand $Q$ as a linear combination of monomials $x^{\alpha}$:
$$
Q=\sum_{\alpha\in {\N}^n,c_{\alpha}\neq 0	} c_{\alpha}x^{\alpha},
$$
where $\alpha=(\alpha_1, \ldots, \alpha_n)\in {\N}^n$, then  it follows from the definition of $Q$ that $0\leq \alpha_i\leq |N|=p-t$ for each $i$.

On the other hand $Q(\ve 0)=\prod_{1\leq i\leq n} \prod_{\alpha\in N} (-\alpha)\neq 0$,
because $0\notin N$.  But it follows from the inclusion (\ref{indir}) that $Q(\ve a_1-\ve a_2)=0$ for each $\ve a_1,\ve a_2\in A$, $\ve a_1\neq \ve a_2$. 

Namely if $\ve a_1,\ve a_2\in A$, $\ve a_1\neq \ve a_2$, then it follows from the inclusion (\ref{indir}) that $\ve a_1-\ve a_2\in ({\Fp})^n\setminus K^n$ and consequently there exists an index $1\leq i\leq n$ such that $(\ve a_1-\ve a_2)_i\notin K$. Hence $(\ve a_1-\ve a_2)_i\in N$ and the definition of $Q$ implies that   $Q(\ve a_1-\ve a_2)=0$.

Consider the polynomials
$$
P_{\ve a}(\ve x):=Q(\ve a-\ve x)\in \Fp[\ve x]
$$
for each $\ve a\in A$. Then it follows from  Proposition \ref{det} that
 $\{P_{\ve a}:~ \ve a\in A\}$ are linearly independent polynomials. Namely the  matrix $B:=(P_{\ve a}(\ve b))_{\ve a,\ve b\in A}$ is a diagonal matrix, where each diagonal entry is nonzero.

On the other hand, if we expand $P_{\ve a}$ as a linear combination of monomials, then  all monomials appearing in this linear combination  contained in the set of monomials 
$$   
\{x_1^{\alpha_1}\cdot\ldots \cdot x_n^{\alpha_n}:~ 0\leq \alpha_i \leq |N| \mbox{ for each } i\}.
$$


Consequently
$$
|A|\leq (|N|+1)^n=(p-t + 1)^n,
$$
a contradiction.
\qed


\end{document}